\documentclass[10pt,twoside,reqno]{amsart}
\usepackage{mathrsfs}
\calclayout

\numberwithin{equation}{section}
\makeatletter

\renewcommand{\@secnumfont}{\bfseries}

\renewcommand{\section}{\@startsection{section}{1}%
  {0mm}{.7\linespacing\@plus\linespacing}{.5\linespacing}
  {\normalfont\bfseries\centering}}

\newcommand{\bibsection}{\@startsection{section}{1}%
  {0mm}{.7\linespacing\@plus\linespacing}{.5\linespacing}
  {\normalfont\scshape\centering}}

\renewcommand{\@biblabel}[1]{#1.}

\newtheorem{thm}{\bf Theorem}[section]
\newtheorem{lem}[thm]{\bf Lemma}

\newtheorem{prop}[thm]{\bf Proposition}

\begin{document}

\vspace{1.3cm}

\title
        {Identities of symmetry for $(h,q)-$extension of higher-order Euler polynomials}

\author{Dae San Kim, Taekyun Kim, Jong Jin Seo }

\thanks{\scriptsize }
\address{1\\ Department of Mathematics\\
             Sogang University\\
              Seoul 121-742, Republic of Korea}
\email{dskim@sogang.ac.kr}
\address{2\\Department of  Mathematics\\
            Kwangwoon University\\
           Seoul 139-701, Republic of Korea}
\email{tkkim@kw.ac.kr}
\address{3\\Department of Applied Mathematics\\
            Pukyong National University\\
            Busan 608-737, Republic of Korea}
\email{seo2011@pknu.ac.kr}

\keywords{multiple $q-$Euler zeta function, $(h,q)-$extension of higher-order Euler polynomials}

\maketitle

\begin{abstract} In this paper, we study some symmetric properties of the multiple $q-$Euler zeta function.
From these properties, we derive several identities of symmetry for the $(h,q)-$extension of higher-order Euler polynomials,
which is an answer to a part of open question in $[7]$.
\end{abstract}

\pagestyle{myheadings}
\markboth{\centerline{\scriptsize Dae San Kim, Taekyun Kim, Jong Jin Seo}}
          {\centerline{\scriptsize Identities of symmetry for $(h,q)-$extension of higher-order Euler polynomials}}
\bigskip
\bigskip
\medskip
\section{\bf Introduction}
\bigskip
\bigskip

Let $\mathbb C$ be the complex number field. We assume that $q\in \mathbb C$ with $|q|<1$ and the $q-$number is defined by
$[x]_q=\frac{1-q^x}{1-q}$. Note that $\lim_{q\rightarrow 1} {[x]_q=x}.$ As is well known, the higher-order Euler polynomials $E_n^{(r)}(x)$ are defined by the generating function to be
\begin{equation}
 F^{(r)} (x,t)=\left(\frac{2}{e^t+1}\right)^r e^{xt}=\sum_{n=0}^\infty E_n^{(r)}(x)\frac{t^n}{n!}, \ \ \textnormal{(see [4], [16]),}
\end{equation}
where $|t|<\pi $.\\
When $x=0, E_n^{(r)}=E_n^{(r)}(0)$ are called the Euler numbers of order $r.$
Recently, the second author defined the $(h,q)-$extension of higher-order Euler polynomials, which is given by the generating function to be
\begin{equation}\begin{split}
 F_q^{(h,r)} (x,t)= &[2]_q^r \sum_{m_1,\cdot\cdot\cdot,m_r=0}^\infty q^{\sum_{j=1}^{r}(h-j+1) m_j} (-1)^{\sum_{j=1}^{r} m_j} e^{[m_1+\cdot\cdot\cdot+m_r+x]_q t}\\
=&\sum_{n=0}^{\infty}E_{n,q}^{(h,r)}(x)\frac{t^n}{n!}, \ \ \textnormal{(see [6], [8]),}
\end{split}\end{equation}
where $h \in {\mathbb Z}$ and $r \in {\mathbb Z}_{\geq 0}.$\\
Note that $\lim_{q\rightarrow 1}F_{q}^{(h,r)}(x,t)=(\frac{2}{e^t+1})^{r}e^{xt}=\sum_{n=0}^{\infty}E_n^{(r)}(x){\frac{t^n}{n!}}.$\\
By $(1.2)$, we get
\begin{equation}
\ F_q^{(h,r)}(t,x)=[2]_{q}^{r}\sum_{m=0}^{\infty} {{m+r-1} \choose {m}}_q(-q^{h-r+1})^m e^{[m+x]_{q}t},
 \ \ \textnormal{(see [6], [8]),}
\end{equation}
 where ${{x} \choose {m}}_{q}= {\frac{[x]_{q} [x-1]_{q} [x-2]_{q} \cdot\cdot\cdot [x-m+1]_{q}}{[m]_{q}!}}.$\\
From $(1.3)$, we can derive the following equation:
\begin{equation}\begin{split}
E_{n,q}^{(h,r)}(x)=&{\frac{[2]_q^r}{(1-q)^n}} \sum_{l=0}^n {n \choose l}{\frac{(-q^x)^l}{(-q^{h-r+l+1}:q)}_r}\\
=&[2]_q^r \sum_{m=0}^{\infty}{m+r-1 \choose m}_q(-q^{h-r+1})^m [m+x]_q^n, \ \ \textnormal{(see [6]),}
\end{split}\end{equation}
 where $(x:q)_n=(1-x)(1-xq)\cdot\cdot\cdot(1-xq^{n-1}).$\\
In $[6]$ and $[8]$, the second author constructed the multiple $q-$Euler zeta function which interpolates the $(h,q)-$ extension of higher-order Euler
polynomials at negative integers as follows :\\
\begin{equation}\begin{split}
\zeta_{q,r}^{(h)}(s,x)=&{\frac{1}{\Gamma (s)}}\int_{0}^{\infty}F_{q}^{(h,r)}(x,t)t^{s-1}dt\\
=&[2]_q^r\sum_{m_1,\cdot\cdot\cdot\ ,m_r=0}^\infty\frac{(-1)^{m_1+\cdot\cdot\cdot+m_r}q^{\sum_{j=1}^{r}(h-j+1)m_j}}{[m_1+\cdot\cdot\cdot+m_r+x]_q^s}, \ \ \textnormal{(see [6]),}
\end{split}\end{equation}
 where ${h,s}\in{\mathbb C}, x\in{\mathbb R}$ with $x\neq{0,-1,-2,\cdot\cdot\cdot}.$\\
From $(1.5)$, we have\\
\begin{equation}\begin{split}
\zeta_{q,r}^{(h)}(s,x)=[2]_q^r\sum_{m=0}^{\infty} {{m+r-1}\choose {m}}_q {(-q^{h-j+1})^m}\frac{1}{[m+x]_q^s} .
\end{split}\end{equation}
Using the Cauchy residue theorem and Laureut series in $(1.5)$, we obtain the following lemma.
\begin{lem}\label{Lemma 1} For $n\in{\mathbb Z}_{\geq 0}$ and $h\in{\mathbb Z},$ we have\\
\begin{equation*}
 \zeta_{q,r}^{(h)}(-n,x)=E_{n,q} ^{(h,r)}(x), \ \ \textnormal{(see [6], [8])}.
\end{equation*}
\end{lem}

In $[7]$, the second author introduced many identities of symmetry for Euler and Bernoulli polynomials which are derived from the $p$-adic integral expression of the generating function and
suggested an open problem about finding identities of symmetry for the Carlitz's type $q$-Euler numbers and polynomials.\\
When $x=0$, $E_{n,q} ^{(h,r)}=E_{n,q} ^{(h,r)}(0)$ are called the $(h,q)$-Euler numbers of order $r$.\\
 From $(1.3)$ and $(1.4)$, we can derive the following equation :\\
\begin{equation}\begin{split}
 E_{n,q} ^{(h,r)}(x)=(q^x E_{q} ^{(h,r)}+[x]_q)^n=\sum_{l=0}^{n} {{n}\choose {l}}q^{lx}  E_{l,q} ^{(h,r)}[x]_q^{n-l} ,
\end{split}\end{equation}
with the usual convention about replacing $(E_q^{(h,r)})^n$ by $E_{n,q} ^{(h,r)}.$\\
Recently, Y. Simsek introduced recurrence symmetric identities for $(h,q)$-Euler polynomials and alternating sums of powers of consecutive $(h,q)$- integers (see $[16]$).\\
In this paper, we investigate some symmetric properties of the multiple $q$-Euler zeta function.
From our investigation, we give some new identities of symmetry for the $(h,q)$-extension of higher-order Euler polynomials, which is an answer to a part of open question in $[7]$.\\

\section{\bf Identities for $(h,q)-$extension of higher-order Euler Polynomials }
\bigskip
\medskip
In this section, we assume that $h\in{\mathbb Z}$ and $a, b\in{\mathbb N}$ with $a\equiv 1(mod\ 2)$ and $b\equiv 1(mod\ 2).$
Now, we observe that\\
\begin{equation}\begin{split}
&\frac{1}{[2]_{q^a}^r}\zeta_{q^a,r}^{(h)} \left(s, bx+\frac{b({j_1+\cdot\cdot\cdot+j_r})}{a}\right)\\
& \ \ \ \ \ =\sum_{m_1,\cdot\cdot\cdot,m_r=0}^\infty \frac{(-1)^{m_1+\cdot\cdot\cdot+m_r}q^{a\sum_{j=1}^{r}(h-j+1)m_j}}
{[m_1+\cdot\cdot\cdot+m_r+bx +\frac{b(j_1+\cdot\cdot\cdot+j_r)}{a}]_{q^a}^s}\\
& \ \ \ \ \ =[a]_q^s\sum_{m_1,\cdot\cdot\cdot,m_r=0}^\infty\frac{q^{a\sum_{j=1}^{r}(h-j+1)m_j}{(-1)^{m_1+\cdot\cdot\cdot+m_r}}}{[b\sum_{l=1}^{r} j_l+{abx}+{a}\sum_{l=1}^{r}{m_l}]_q^s}\\
& \ \ \ \ \ =[a]_q^s\sum_{m_1,\cdot\cdot\cdot,m_r=0}^\infty \sum_{i_1,\cdot\cdot\cdot,i_r=0}^{b-1}\frac{(-1)^{\sum_{l=1}^{r} (i_l+bm_l)} q^{a\sum_{j=1}^{r} (h-j+1)(i_j+m_jb)}}{[ab(x+\sum_{l=1}^{r} m_l) +b\sum_{l=1}^{r} j_l +a\sum_{l=1}^{r} i_l]_q^s}.\\
\end{split}\end{equation}
Thus, by $(2.1)$, we get
\begin{equation}\begin{split}
&\frac{[b]_q^s}{[2]_{q^a}^r} \sum_{j_1,\cdot\cdot\cdot,j_r=0}^{a-1} (-1)^{\sum_{l=1}^{r}j_l} q^{b\sum_{l=1}^{r}(h-l+1)j_l} \zeta_{{q^a},r}^{(h)} \left(s,bx+\frac{b(j_1+\cdot\cdot\cdot\,+j_r)}{a}\right)\\
=&[a]_q^s [b]_q^s \sum_{j_1,\cdot\cdot\cdot,j_r=0}^{a-1} \sum_{i_1,\cdot\cdot\cdot,i_r=0}^{b-1} \sum_{m_1,\cdot\cdot\cdot,m_r=0}^\infty \frac{(-1)^{\sum_{l=1}^{r} (i_l+j_l+m_l)} q^{a\sum_{l=1}^{r}(h-l+1)i_l+b\sum_{l=1}^{r}(h-l+1)j_l}}{[ab(x+\sum_{l=1}^{r}m_l)+\sum_{l=1}^{r}(bj_l+ai_l)]_q^s}\\
&\times q^{ab\sum_{l=1}^{r} m_l}.
\end{split}\end{equation}
By the same method as $(2.2)$, we see that
\begin{equation}\begin{split}
&\frac{[a]_q^s}{[2]_{q^b}^r} \sum_{j_1,\cdot\cdot\cdot,j_r=0}^{b-1} (-1)^{\sum_{l=1}^{r}j_l} q^{a\sum_{l=1}^{r}(h-l+1)j_l} \zeta_{{q^b},r}^{(h)} \left(s,ax+\frac{a(j_1+\cdot\cdot\cdot\,+j_r)}{b}\right)\\
=&[b]_q^s [a]_q^s \sum_{j_1,\cdot\cdot\cdot,j_r=0}^{b-1} \sum_{i_1,\cdot\cdot\cdot,i_r=0}^{a-1} \sum_{m_1,\cdot\cdot\cdot,m_r=0}^\infty \frac{(-1)^{\sum_{l=1}^{r} (i_l+j_l+m_l)} q^{a\sum_{l=1}^{r}(h-l+1)j_l+b\sum_{l=1}^{r}(h-l+1)i_l}}{[ab(x+\sum_{l=1}^{r}m_l)+\sum_{l=1}^{r}(bi_l+aj_l)]_q^s}\\
&\times q^{ab\sum_{l=1}^{r} m_l}.
\end{split}\end{equation}
Therefore, by$(2.2)$ and $(2.3)$, we obtain the following theorem.
\bigskip
\begin{thm}\label{Theorem 1.} For $a,b \in{\mathbb N}$, with $a\equiv1(mod\ 2)$ and $b\equiv1(mod\ 2)$, we have\\
\begin{equation*}\begin{split}
&[2]_{q^b}^r [b]_q^s \sum_{j_1,\cdot\cdot\cdot,j_r=0}^{a-1} (-1)^{j_1+\cdot\cdot\cdot+j_r} q^{b\sum_{l=1}^{r} (h-l+1)j_l} \zeta_{{q^a},r}^{(h)} \left(s,bx+\frac{b(j_1+\cdot\cdot\cdot+j_r)}{a}\right)\\
&=[2]_{q^a}^r [a]_q^s \sum_{j_1,\cdot\cdot\cdot,j_r=0}^{b-1} (-1)^{j_1+\cdot\cdot\cdot+j_r} q^{a\sum_{l=1}^{r} (h-l+1)j_l} \zeta_{{q^b},r}^{(h)} \left(s,ax+\frac{a(j_1+\cdot\cdot\cdot+j_r)}{b}\right).\\
\end{split}\end{equation*}
\end{thm}
\bigskip
From Lemma $1.1$ and Theorem $2.1$, we can derive the following theorem.

\begin{thm}\label{Theorem 2.} For $n \in{\mathbb Z_{\geq0}}$ and $a,b\in{\mathbb N}$, with $ a\equiv1 (mod\ 2)$ and $b\equiv1 (mod\ 2)$, we have\\
\begin{equation*}\begin{split}
&[2]_{q^b}^r [a]_q^n \sum_{j_1,\cdot\cdot\cdot,j_r=0}^{a-1} (-1)^{j_1+\cdot\cdot\cdot+j_r} q^{b\sum_{l=1}^{r} (h-l+1)j_l} E_{n,{q^a}}^{(h,r)}
 \left(bx+\frac{b(j_1+\cdot\cdot\cdot+j_r)}{a} \right)\\
&=[2]_{q^a}^r [b]_q^n \sum_{j_1,\cdot\cdot\cdot,j_r=0}^{b-1} (-1)^{j_1+\cdot\cdot\cdot+j_r} q^{a\sum_{l=1}^{r} (h-l+1)j_l} E_{n,{q^b}}^{(h,r)} \left(ax+\frac{a(j_1+\cdot\cdot\cdot+j_r)}{b}\right).\\
\end{split}\end{equation*}
\end{thm}

By $(1.4)$, we easily see that
\begin{equation}\begin{split}
E_{n,q}^{(h,k)}(x+y)&=(q^{x+y}E_q^{(h,k)} +[x+y]_q)^n =(q^{x+y}E_q^{(h,k)}+q^{x}[y]_q+[x]_q)^n\\
&=\left(q^x(q^yE_q^{(h,k)}+[y]_q)+[x]_q\right)^n=\sum_{i=0}^{n} {{n}\choose {i}} q^{ix} E_{i,q}^{(h,k)}(y)[x]_q^{n-i}.
\end{split}\end{equation}
Therefore, by $(2.4)$, we obtain the following proposition.

\begin{prop}\label{Proposition 3.} For $n\geq0$, we have\\
\begin{equation*}\begin{split}
E_{n,q}^{(h,k)}(x+y)&=\sum_{i=0}^{n} {n\choose i} q^{ix}E_{i,q}^{(h,k)}(y)[x]_{q}^{n-i}=\sum_{i=0}^{n} {{n}\choose {i}} q^{(n-i)x}E_{n-i,q}^{(h,k)}(y)[x]_{q}^{i}.\\
\end{split}\end{equation*}
\end{prop}

From Proposition$2.3$, we note that
\begin{equation}\begin{split}
&\sum_{j_1,\cdot\cdot\cdot,j_r=0}^{a-1} (-1)^{j_1+\cdot\cdot\cdot+j_r} q^{b\sum_{l=1}^{r} (h-l+1)j_l} E_{n,{q^a}}^{(h,r)}\left(bx+\frac{b(j_1+\cdot\cdot\cdot+j_r)}{a}\right)\\
&=\sum_{j_1,\cdot\cdot\cdot,j_r=0}^{a-1} (-1)^{j_1+\cdot\cdot\cdot+j_r} q^{b\sum_{l=1}^{r} (h-l+1)j_l} \sum_{i=0}^{n} {n\choose i} q^{ia\left(\frac{b(j_1+\cdot\cdot\cdot+j_r)}{a}\right)}E_{i,{q^a}}^{(h,r)}(bx)\\
&\times \left[\frac{b(j_1+\cdot\cdot\cdot+j_r)}{a}\right]_{q^a}^{n-i}\\
&=\sum_{j_1,\cdot\cdot\cdot,j_r=0}^{a-1} (-1)^{j_1+\cdot\cdot\cdot+j_r} q^{b\sum_{l=1}^{r} (h-l+1)j_l} \sum_{i=0}^{n} {n\choose i} q^{(n-i){b(j_1+\cdot\cdot\cdot+j_r)}}E_{n-i,{q^a}}^{(h,r)}(bx)\\
&\times \left[\frac{b(j_1+\cdot\cdot\cdot+j_r)}{a}\right]_{q^a}^{i}\\
&=\sum_{i=0}^{n} {n\choose i} \left(\frac{[b]_q}{[a]_q}\right)^i E_{n-i,{q^a}}^{(h,r)}(bx)\sum_{j_1,\cdot\cdot\cdot\,j_r=0}^{a-1} (-1)^{j_1+\cdot\cdot\cdot+j_r} q^{b\sum_{l=1}^{r} (h+n-l-i+1)j_l}\left[j_1+\cdot\cdot\cdot+j_r\right]_{q^b}^i\\
&=\sum_{i=0}^{n} {n\choose i} \left(\frac{[b]_q}{[a]_q}\right)^i E_{n-i,{q^a}}^{(h,r)}(bx) S_{n,i,{q^b}}^{(h,r)}(a),\\
\end{split}\end{equation}
\begin{equation}\begin{split} where\ \
S_{n,i,q}^{(h,r)}(a)=\sum_{j_1,\cdot\cdot\cdot\,j_r=0}^{a-1} (-1)^{j_1+\cdot\cdot\cdot+j_r} q^{\sum_{l=1}^{r} (h+n-l-i+1)j_l}\left[j_1+\cdot\cdot\cdot+j_r\right]_{q}^{i}.\\
\end{split}\end{equation}
By $(2.5)$, we get
\begin{equation}\begin{split}
&[2]_{q^b}^r [a]_q^n \sum_{j_1,\cdot\cdot\cdot,j_r=0}^{a-1} (-1)^{j_1+\cdot\cdot\cdot+j_r} q^{b\sum_{l=1}^{r} (h-l+1)j_l} E_{n,{q^a}}^{(h.r)}\left(bx+\frac{b(j_1+\cdot\cdot\cdot+j_r)}{a}\right)\\
&=[2]_{q^b}^r \sum_{i=0}^{n} {n\choose i} [a]_q^{n-i}[b]_q^i E_{n-i,{q^a}}^{(h,r)}(bx)S_{n,i,{q^b}}^{(h,r)}(a).\\
\end{split}\end{equation}
By the same method as $(2.7)$, we see that
\begin{equation}\begin{split}
&[2]_{q^a}^r [b]_q^n \sum_{j_1,\cdot\cdot\cdot\,j_r=0}^{b-1} (-1)^{j_1+\cdot\cdot\cdot+j_l} q^{a\sum_{l=1}^{r} (h-l+1)j_l} E_{n,{q^b}}^{(h,r)}\left(ax+\frac{a(j_1+\cdot\cdot\cdot+j_r)}{b}\right)\\
&=[2]_{q^a}^r \sum_{i=0}^{n} {n\choose i} [b]_q^{n-i}[a]_q^i E_{n-i,{q^b}}^{(h,r)}(ax)S_{n,i,{q^a}}^{(h,r)}(b).\\
\end{split}\end{equation}
Therefore, by $(2.7)$ and $(2.8)$, we obtain the following theorem.\\
\begin{thm}\label{Theorem 2.4.} For ${a,b}\in\mathbb{N}$ with $a\equiv 1(mod\ 2)$  and $b\equiv 1(mod\ 2), n\in{\mathbb Z}_{\geq 0}$,\\
\\ let
\begin{equation*}
 S_{n,i,q}^{(h,r)}(a)={\sum_{j_1,\cdot\cdot\cdot,j_r=0}^{a-1}} (-1)^{j_1+\cdot\cdot\cdot+j_r}q^{{\sum_ {l=1}^{r}}(h+n-l-i+1)j_l}[j_1+\cdot\cdot\cdot+j_r]_q^i.
\end{equation*} Then we have
\begin{equation*}
[2]_{q^b}^r \sum_{i=0}^{n}{n\choose i}[a]_q^{n-i}[b]_q^iE_{n-i,q^a}^{(h,r)}(bx)S_{n,i,q^b}^{(h,r)}(a)
=[2]_{q^a}^r \sum_{i=0}^{n}{n\choose i}[b]_q^{n-i}[a]_q^iE_{n-i,q^b}^{(h,r)}(ax)S_{n,i,q^a}^{(h,r)}(b).
\end{equation*}
\end{thm}
It is not difficult to show that
\begin{equation}
[x+y+m]_q(u+v)-[x]_qv=[x]_qu+q^x[y+m]_q(u+v).
\end{equation}
From $(2.9)$, we note that\\
\begin{equation}\begin{split}
&e^{[x]_qu}{\sum_{m_1,\cdot\cdot\cdot,m_r=0}^\infty} q^{{\sum_ {j=1}^{r}}(h-j+1)m_j}(-1)^{{\sum_ {j=1}^{r}}m_j} e^{{[m_1+\cdot\cdot\cdot+m_r+y]_q}q^x(u+v)}\\
&= e^{-[x]_qv}{\sum_{m_1,\cdot\cdot\cdot,m_r=0}^\infty} q^{{\sum_ {j=1}^{r}}(h-j+1)m_j}(-1)^{{\sum_ {j=1}^{r}}m_j} e^{{[x+y+m_1+\cdot\cdot\cdot+m_r]_q}(u+v)}.
\end{split}\end{equation}
The left hand side of $(2.10)$ multiplied by $[2]_q^r$ is given by\\
\begin{equation}\begin{split}
&[2]_q^re^{[x]_qu}{\sum_{m_1,\cdot\cdot\cdot,m_r=0}^\infty} q^{{\sum_ {j=1}^{r}}(h-j+1)m_j}(-1)^{{\sum_ {j=1}^{r}}m_j} e^{{[m_1+\cdot\cdot\cdot+m_r+y]_q}q^x(u+v)}\\
& \ \ \ \ = e^{[x]_qu}{\sum_{n=0}^\infty} q^{nx} E_{n,q}^{(h,r)}(y)\frac{1}{n!}(u+v)^n\\
& \ \ \ \ =\left({\sum_{l=0}^\infty} [x]_q^l\frac{u^l}{l!}\right) \left({\sum_{n=0}^\infty}q^{nx} E_{n,q}^{(h,r)}(y) {\sum_{k=0}^n}\frac{u^k}{k!(n-k)!}v^{n-k}\right)\\
& \ \ \ \ =\left({\sum_{l=0}^\infty} [x]_q^l\frac{u^l}{l!}\right) \left({\sum_{k=0}^\infty}{\sum_{n=0}^\infty}q^{(n+k)x} E_{n+k,q}^{(h,r)}(y) \frac{u^k}{k!}\frac{v^n}{n!}\right)\\
& \ \ \ \ =\sum_{m=0}^\infty \sum_{n=0}^\infty \left({\sum_{k=0}^m}{m\choose k}q^{(n+k)x}  E_{n+k,q}^{(h,r)}(y)[x]_q^{m-k}\right)  \frac{u^m}{m!}\frac{v^n}{n!}
\end{split}\end{equation}
The right hand side of $(2.10)$ multiplied by $[2]_q^r$ is given by\\
\begin{equation}\begin{split}
&[2]_q^re^{-[x]_qv}{\sum_{m_1,\cdot\cdot\cdot,m_r=0}^\infty} q^{{\sum_ {j=1}^{r}}(h-j+1)m_j}(-1)^{{\sum_ {j=1}^{r}}m_j} e^{{[x+y+m_1+\cdot\cdot\cdot+m_r]_q}(u+v)}\\
& \ \ \ \ = e^{-[x]_qv}{\sum_{n=0}^\infty}  E_{n,q}^{(h,r)}(x+y)\frac{1}{n!}(u+v)^n\\
& \ \ \ \ =\left({\sum_{l=0}^\infty} \frac{(-[x]_q)^l}{l!}v^l\right)\left({\sum_{m=0}^\infty}{\sum_{k=0}^\infty} E_{m+k,q}^{(h,r)}(x+y) \frac{u^m}{m!}\frac{v^k}{k!}\right)\\
& \ \ \ \ =\sum_{n=0}^\infty \sum_{m=0}^\infty \left({\sum_{k=0}^n}{n\choose k}  E_{m+k,q}^{(h,r)}(x+y)(-[x]_q)^{n-k}\right)  \frac{u^m}{m!}\frac{v^n}{n!}\\
& \ \ \ \ =\sum_{n=0}^\infty \sum_{m=0}^\infty \left({\sum_{k=0}^n}{n\choose k}  E_{m+k,q}^{(h,r)}(x+y)q^{(n-k)x}[-x]_q^{n-k}\right)  \frac{u^m}{m!}\frac{v^n}{n!}.
\end{split}\end{equation}
Therefore, by $(2.10)$, $(2.11)$  and $(2.12)$, we obtain the following theorem.\\
\begin{thm}\label{Theorem 2.5.} For $m,n\geq 0$, we have
\begin{equation*}
\sum_{k=0}^m {m \choose k} q^{(n+k)x}E_{n+k,q}^{(h,r)}(y)[x]_q^{m-k}
=\sum_{k=0}^n {n \choose k} E_{m+k,q}^{(h,r)}(x+y)q^{(n-k)x}[-x]_q^{n-k}.
\end{equation*}
\end{thm}
\bigskip
Remark. Recently, several authors have studied $(h,q)-$extension of Bernoulli and Euler polynomials $\textnormal{(see[1]-[5],\ [9]-[17])}.$

\end{document}